\documentclass[letterpaper,twoside,12pt,leqno]{article}
\usepackage{amsmath,amsthm,amssymb, enumerate}
            \newcommand{\marginalnote}[1]{}
\usepackage[all]{xy}
\usepackage{graphicx}
\usepackage{epsfig}
\usepackage{color}	

\swapnumbers \theoremstyle{plain}

\newtheorem{theorem}{Theorem}[section]

\newtheorem{lemma}[theorem]{Lemma}
\newtheorem{corollary}[theorem]{Corollary}
\theoremstyle{remark}

\theoremstyle{definition}
\newtheorem{remark}[theorem]{Remark}
\newtheorem{example}[theorem]{Example}
\newtheorem{Rems}[theorem]{Remarks}
\newtheorem{Rev}[theorem]{Review}

\newtheorem{definition}[theorem]{Definition}
\newtheorem{definitions}[theorem]{Definitions}

\newtheorem{observation}[theorem]{Observation}

\newcommand{\rconj}[2]{{}^{#1}\mkern-1mu{#2}}
\newcommand{\gen}[1]{\left\langle#1\right\rangle}

\newcommand{\gp}[2]{\gen{#1\mid #2}}
\newcommand{\prs}[2]{\gen{#1\parallel #2}}
\newcommand{\ol}{\overline}
\newcommand{\Z}{\mathbb{Z}}
\newcommand{\R}{\mathbb{R}}

\newcommand{\ball}{\mathbb B}

\def\coloneqq{\mathrel{\mathop\mathchar"303A}\mkern-1.2mu=}

\DeclareMathOperator{\col}{color}
\DeclareMathOperator{\dist}{d}
\DeclareMathOperator{\diam}{diam} 
\DeclareMathOperator{\adim}{adim} 
\DeclareMathOperator{\Star}{star}
\begin{document}

\pagestyle{myheadings} \markboth{}{Yago Antol\'{i}n Pichel}
\title{On Cayley graphs of virtually free groups.}
\date{\today}
\author{Yago Antol\'{i}n Pichel}
\maketitle
\begin{abstract}
In 1985, Dunwoody showed that finitely presentable groups are accessible. Dunwoody's result was used to show that context-free groups, groups quasi-isometric to trees or finitely presentable groups of asymptotic dimension 1 are virtually free.

Using another theorem of Dunwoody of 1979, we study when a group is virtually free in terms of its Cayley graph, and we obtain new proofs of the mentioned results and others previously depending on these.

\medskip

{\footnotesize
\noindent \emph{2000 Mathematics Subject Classification.} Primary: 20F65 ;
Secondary: 05C25
\noindent \emph{Key words.} Virtually free group, accessibility, Cayley graph, asymptotic dimension.}
\end{abstract}
\maketitle

\section{Notation and background}
Let $G$ be a discrete, finitely generable, multiplicative group, fixed throughout the article.

Let $\Gamma= (\Gamma,V\Gamma, E\Gamma, \iota, \tau)$ be an oriented graph fixed throughout the article.

Recall that $\Gamma$ is a set with a specified non-empty subset $V\Gamma,$ $E\Gamma$ is the complement of $V\Gamma$ in $\Gamma$, and $\iota,\tau$ are functions from $E\Gamma$ to $V\Gamma.$  The elements of $V\Gamma$ are called vertices, the elements of $E\Gamma$ are called edges and the functions $\iota,\tau\colon E\Gamma\to V\Gamma$ are called incidence functions.%

For two subsets $A,$ $B$ of a set $S,$ the complement of $A \cap B$ in $A$ will be denoted by $A-B$ (and not by $A\setminus B$ since we let $G\backslash Y$ denote the set of $G$-orbits of a left $G$-set $Y$).

A  \textit{sequence} is a \textit{set} endowed with a specified listing of its elements,
usually represented as a vector in which the coordinates are the elements of
(the underlying set of) the sequence.
For two sequences $A,$ $B,$  their  concatenation will be denoted~$A \vee B.$

We will find it useful to have a notation  for intervals in $\Z$ that is different from the notation for intervals in $\mathbb{R}.$

Let $i,$ $j,$ $n \in \Z.$

We write \begin{equation*}[i{\uparrow}j]\coloneqq
 \begin{cases}
(i,i+1,\ldots, j-1, j) \in \Z^{j-i+1}   &\text{if $i \le j,$}\\
() \in \Z^0 &\text{if $i > j.$}
\end{cases}
\end{equation*}
Also, $[i{\uparrow}\infty[\,\,\,\, \coloneqq (i,i+1,i+2,\ldots)$ and
$[i{\uparrow}\infty]\coloneqq  [i{\uparrow}\infty[ \,\,\,\vee\,\, \{\infty\}.$
We define $[j{\downarrow}i]$ to be the reverse of the sequence $[i{\uparrow}j],$
that is, $(j,j-1,\ldots,i+1,i).$

We shall use sequence notation to define families of
indexed symbols.
Let $v$ be a symbol.
For each $k \in \Z,$ we let $v_k$ denote the ordered pair
$(v,k).$  We let \begin{equation*}v_{[i{\uparrow}j]}\coloneqq  \begin{cases}
(v_i,v_{i+1},   \cdots, v_{j-1},   v_j)  &\text{if $i \le j,$}\\
()  &\text{if $i > j.$}
\end{cases}\end{equation*}
Also, $v_{[i{\uparrow}\infty[\,} \coloneqq (v_i,v_{i+1},v_{i+2},\ldots).$
We define $v_{[j{\downarrow}i]}$ to be the reverse of the sequence~$v_{[i{\uparrow}j]}.$

Now suppose that $v_{[i{\uparrow}j]}$ is a sequence \textit{in} the group $G,$
that is, there is a specified map $v_{[i{\uparrow}j]} \to G.$  We treat the
elements of $v_{[i{\uparrow}j]}$ as elements of $G,$ possibly with repetitions,
and we define
\begin{align*}\Pi v_{[i{\uparrow}j]}&\coloneqq  \begin{cases}
v_i   v_{i+1}   \cdots v_{j-1}   v_j \in G   &\text{if $i \le j,$}\\
1 \in G &\text{if $i > j.$}
\end{cases}\\ \Pi v_{[j{\downarrow}i]} &\coloneqq  \begin{cases}
v_j   v_{j-1}   \cdots v_{i+1}   v_i \in G   &\text{if $j \ge i,$}\\
1 \in G &\text{if $j< i .$}
\end{cases}\end{align*}

Let $S$ be a set. Suppose that $v_{[i\uparrow j]}$ is a sequence in the set of subsets of $S,$ then we define $\cap v_{[i\uparrow j]}\coloneqq \cap_{k\in [i\uparrow j]}v_k$ and $\cup v_{[i\uparrow j]}\coloneqq \cup_{k\in [i\uparrow j]}v_k.$

Let $E\Gamma^{\pm 1}\coloneqq E\Gamma\times \{-1,1\}$ and for $e\in E\Gamma,$ let $e^{-1}$ denote $(e,-1),$ $e=e^{1}$ denote $(e,1)$ and we extend $\iota,\tau$ to $E\Gamma^{\pm 1}$ defining $\iota e^{-1}=\tau e$ and $\tau e^{-1}=\iota e.$

A {\it path } $\gamma$  {\it in} $\Gamma$ is a pair of sequences $(v_{[0\uparrow n]},e_{[1\uparrow n]})$, where $v_{[0\uparrow n]}$ is a sequence in $V\Gamma,$  $e_{[1\uparrow n]}$ is a sequence in $E\Gamma^{\pm 1}$ such that for $i\in [1\uparrow n],$ $\iota e_i=v_{i-1}$ and $\tau e_i= v_{i}.$ The {\it length} of $\gamma$  is $n.$ The path is {\it closed} if $v_0=v_n.$ The path $\gamma$ is {\it geodesic} if for any other path $(u_{[0\uparrow m]},f_{[1\uparrow m]})$ in $\Gamma$ with $u_0=v_0$ and $v_n=u_m,$ $n\leq m.$

A sequence $v_{[0\uparrow n]}$ in $V\Gamma,$ is {\it a  path in $V\Gamma$} if there exists a sequence $e_{[1\uparrow n]}$ in $E\Gamma^{\pm 1}$ such that  $(v_{[0\uparrow n]},e_{[1\uparrow n]})$ is a path in $\Gamma.$ 

We use $\dist_\Gamma$ or simply $\dist$ to denote the path metric on $V\Gamma.$ That is, $\dist_\Gamma(u,v)= n$ if there exists a geodesic path  $(v_{[0\uparrow n]},e_{[1\uparrow n]})$ in $\Gamma$ with $v_0=u$ and $v_n=v$ and $\dist_\Gamma(u,v)=\infty$ if no such geodesic path exists.

Let $X$ be a finite generating set of $G.$ The graph $\Gamma$ is the \textit{Cayley graph} of $G$ with respect to $X$ if $V\Gamma=G,$ $E\Gamma=G\times X$ and for $(g,x)\in E\Gamma,$ $\iota (g,x)=g$ and $\tau(g,x)=g\cdot x.$ In this case, for $g\in G$, we put $|g|\coloneqq\dist_\Gamma(1,g).$

\section{Introduction and outline}
A group is \textit{virtually free} if it has a free subgroup of finite index. 

By Stalling's ends theorem, $G$ has \textit{more than one end} if and only if there exists a $G$-tree with finite edge stabilizers such that no vertex is stabilized by $G$.

The group $G$ is \textit{accessible} if there exists a $G$-tree with  finite edge stabilizers, and each vertex stabilizer having at most one end.

The group $G$ is \textit{context-free} if for every finite generating set $X$ of $G,$ the set \begin{equation*}L(G,X)=\{x_{[1\uparrow n]}\text{ sequence in }X^{\pm 1}: \Pi x_{[1\uparrow n]}=1\}\end{equation*} is a context-free language. See Review~\ref{Rev:lang} for a complete definition.

In 1983, Muller and Schupp \cite[Theorem III]{MullerSchupp83} showed that $G$ is context-free and accessible if and only if $G$ is virtually free.

By \cite[Theorem I]{MullerSchupp83},  a group is context-free if every closed path in the Cayley graph is $k$-triangulable. In particular it follows that being context-free can be expressed in terms of the Cayley graph and that context-free groups are finitely presentable.

In 1985, Dunwoody \cite{Dunwoody85} proved that finitely presented groups are accessible,
and so the accessibility hypothesis  of \cite[Theorem III]{MullerSchupp83} can be removed.

In this article we investigate how to characterize virtual freeness of a group 
in terms of its Cayley graph. Using an accessiblility  result (Theorem~\ref{T:Dun}) which is implicit in an earlier paper of Dunwoody \cite{Dunwoody79} of 1979 we obtain new proofs of several results which originally used the accessibility of finitely presented groups.  Let us list these results.

Let $\Gamma$ be the Cayley graph of $G,$ then $G$ is virtually free if and only if one of the following hold:

\begin{enumerate}
\item[(1)] There exists $m\geq 0$ such that every closed path in $\Gamma$ is $m$-triangulable (see \cite[Theorems I and III]{MullerSchupp83}).
\item[(2)] $\Gamma$ is quasi-isometric to a Cayley graph of a free group (see \cite[7.19]{GhysHarpe}).
\item[(3)] $G$ is finitely presentable and $\adim \Gamma =1$ (see \cite[Cor 1.2]{FujiwaraWhyte},\cite[Thm 2]{Gentimis} or \cite[Cor G]{JanuszkiewiczSwiatkowski}).
\item[(4)] $\Gamma$ admits a uniform spanning tree (\cite{Woess}).
\item[(5)] The ends of $\Gamma$ have finite diameter (\cite{Woess}).
\item[(6)] $\Gamma$ has finite strong tree-width (\cite{KuskeLohrey}).
\end{enumerate}

The original proofs of (1)-(3) directly invoke the accesibility of finitely presentable groups. The proofs of (4)-(6) invoke (1). As (1) is equivalent to $G$ being context-free, we found  two other characterizations of $G$ being virtually free: 

\begin{enumerate}
\item[(7)] There exists a generating set of $G$ and $k\geq 0$ such that every $k$-locally geodesic in $G$ is a geodesic (\cite{GilmanHermillerHoltRees})
\item[(8)] There exists a generating set of $G$ such that $G$ admits a (regular) broomlike combing with respect to the generating set(\cite{BridsonGilman}).
\end{enumerate}

In section~\ref{S:Dun}, we will review the classical geometric characterization of $G$ being virtually free in terms of actions on trees. We will use an accessibility result of Dunwoody of 1979, implicit in \cite{Dunwoody79}.

In section~\ref{S:graphs}, we show the equivalence of many conditions for graphs. Using results of section~\ref{S:Dun}, we will show that the Cayley graph of $G$  satisfies any of these conditions if and only if $G$ is virtually free. In particular, we give new proofs of (1),(2),(4)-(6).

At the end of the section, we recall the definition of a context-free group in Review~\ref{Rev:lang}, and we make some remarks about the connection of the $m$-triangulations with the characterizations (7) and (8).

In section~\ref{S:adim}, we review the concept of asymptotic dimension. We use (1) to give a proof of (3).

\section{Classical approach}\label{S:Dun}
In this section, using a theorem of Dunwoody, we show that $G$ is virtually free if and only if $G$ acts on a tree with finite stabilizers.

\begin{definitions} Let $S$ be a set and let $(G,S)$ denote the maps from $G$ to $S.$ Let $\alpha\in (G,S)$ and $g\in G.$ There is a natural right $G$-set structure on $(G,S)$ defined by $\alpha g(h)=\alpha(hg)$ for all $h\in G.$

The set $R$ is \textit{almost contained} in $S,$ written $R\subset_a S,$ if $\{r\in R : r\notin S\}$ is finite. If $R\subset_a S$ and $S\subset_a R,$ then the sets $R$ and $S$ are \textit{almost equal}, and it is denoted by $R=_a S.$

Let $S$ be a right $G$-set. A subset $R\subset S$ is  \textit{almost-right-$G$-invariant} if  for all $g\in G,$ $Rg=_a R.$ 

A function $\alpha\in (G,S)$ is \textit{almost-right-$G$-invariant} if for all $s\in S,$ $\alpha^{-1}(s)$ is almost-right-$G$-invariant.
\end{definitions}
Stallings' ends theorem relates almost invariant sets of $G$ with graph of groups decompositions of $G.$ In \cite{KarrasPietrowskiSolitar73}, using Stallings' theorem, Karrass, Pietrowski and Solitar showed that finitely generable virtually free groups act on trees with finite stabilizers. We will follow the approach to the theory of Dunwoody \cite{Dunwoody79}, using the following result.

\begin{theorem}\label{T:Dun} Let $G$ be a finitely generable group and $S$ a finite set. Let $\alpha \colon G\to S$ be almost-right-$G$-stable, then $G$ is the fundamental group of a graph of groups $(G(-),Y)$ where $Y$ is finite, the edge groups are finite and each vertex group is the right-stabilizer of some $\beta=_a\alpha.$
\end{theorem}

This result appears implicit in \cite{Dunwoody79}, as a combination of Lemma 4.7 and Lemma 4.6 with  $H=1.$ There is an explicit formulation in \cite[III.3.4]{Dicks80}. It can also be shown as a consequence of the Almost Stability Theorem \cite[III.8.5]{DicksDunwoody89}.

\begin{theorem}\label{T:A}
Let $G$ be a finitely generable group. The following are equivalent
\begin{itemize}

\item[(A1)] $G$ acts on a tree with finite stabilizers.
\item[(A2)] $G$ is the fundamental group of a finite graph of finite groups.
\item[(A3)] $G$ is virtually free.
\item[(A4)] There exists a finite set $S$ and an almost-right-$G$-invariant function $\alpha\in (G,S)$ such that for all $\beta=_a \alpha$ the right-stabilizer of $\beta$ is finite. 
\end{itemize}
\end{theorem}

\begin{proof}
(A1) and (A2) are equivalent by Bass-Serre theory, see \cite[I.4.1]{DicksDunwoody89} and  \cite[I.7.6]{DicksDunwoody89} for a proof.
(A2) implies (A3) by a classic result of Karrass, Pietrowski and Solitar \cite{KarrasPietrowskiSolitar73}, see \cite[I.7.4]{DicksDunwoody89} for a proof.
(A4) implies (A2) by Theorem~\ref{T:Dun}.

To complete the proof, it is enough to show that (A3) implies (A4). 

Let $g,h\in G$ and suppose first that $G$ is free.

Let $X$ be a free generating set for $G,$  and $S=\{1\}\cup X\cup X^{-1}.$ Let $\alpha\in (G,S)$ be the function that assigns to each $g\in G-\{1\}$ the first generator of the free reduced word in $X$ representing $g$ and $\alpha(1)=1.$

Let $\Gamma$ be the Cayley graph of $G$ with respect to $X.$ Geometrically, $\alpha$ is the function that collapses each connected component of \begin{equation*}\Gamma-\{(1,x),(x,x^{-1}) : x\in X\}\end{equation*} to a vertex.

Suppose that $|g|<|h|.$ Then $(\alpha g)(h)=\alpha(hg)=\alpha(h).$ Hence $\alpha g =_a\alpha.$ 

Let $\beta \in (G,S),$ $\beta=_a\alpha$ and assume $g\neq 1.$  We will show that $g\notin G_{\beta},$ the $G$-stabilizer of $\beta.$

If $g$ is cyclically reduced, $\alpha(g)\neq \alpha({g}^{-1})$ then \begin{equation*}\beta(g^n)=\alpha(g^n)\neq \alpha(g^{-n})=\beta(g^{-n})\end{equation*} for infinitely many $n\in [1\uparrow\infty[,$  and thus $(\beta g^{2n})(1)\neq \beta(1)$ for infinitely many $n$. Hence $g\notin G_\beta.$

If $g$ is not cyclically reduced, then $\alpha(g)=\alpha(g^{-1})$ and  there exist $g_1,k\in G-\{1\}$ such that $g=kg_1{k}^{-1}$ and $\alpha(g_1)\neq\alpha(g_1^{-1}).$ By the previous discussion, $g_1\notin G_{(\beta k)}$ and hence $g\notin G_\beta$.

We consider now the general case. Let $H$ be a finite-index free subgroup of $G,$ $S$ a finite set and $\alpha\in (H,S)$ an almost-right-$H$-invariant function such that for all $\beta=_a\alpha,$ $H_\beta$ is finite.

Let $R$ be a transversal for the left $H$-multiplication in $G.$ For each $g\in G$, denote by $h_g$ and $r_g$ the unique elements $h_g\in H$ and $r_g\in R$ such that $g=h_gr_g$. Let $\widehat{\alpha} \colon G\to S$ be defined by $\widehat{\alpha}(g)=\alpha(h_g).$

For $s\in S,$
\begin{equation*}\widehat{\alpha}^{-1}(s)g=(\bigcup_{r\in R}\alpha^{-1}(s)r)g=\bigcup_{t\in Rg}\alpha^{-1}(s)h_{t} r_{t}=_a \bigcup_{t\in Rg}\alpha(s)^{-1} r_{t}=\widehat{\alpha}^{-1}(s).\end{equation*}

Hence $\widehat{\alpha}$ is almost-right-$G$-stable. 

Let $\widehat{\beta}=_a\widehat{\alpha}.$ For all but a finite number of $g\in G_{\widehat{\beta}}$ \begin{equation*} \widehat{\beta} g=\widehat{\beta} h_g=\widehat{\beta}.\end{equation*} Let $\beta$ be the restriction of $\widehat{\beta}$ to $H$. Then for all but a finite number of $g\in G_{\widehat{\beta}}$ $\beta h_g=\beta.$ As $\beta=_a \alpha$, it follows that $G_{\widehat{\beta}}$ is finite.
\end{proof} 

\section{Equivalent graph conditions}\label{S:graphs}
Throughout $T=(VT,ET,\iota_T,\tau_T)$ will be a fixed graph. 

We show the equivalence of many conditions for $\Gamma$ in a succession of Lemmas. We summarize the results in Theorem~\ref{T:Eqcond}. Furthermore, for Cayley graphs, these conditions characterize virtual freeness.

\begin{definitions}
Let $(X,\dist_X)$ and $(Y,\dist_Y)$ be two metric spaces.
Let $\phi \colon X\to Y$ be a map.

The map $\phi$ is \textit{Lipschitz continuous} if there exists a constant $\lambda>0$ such that for all $x_1,x_2\in X,$ $\dist_Y(\phi(x_1),\phi(x_2))\leq \lambda \dist_X(x_1,x_2).$

The map $\phi$ is \textit{large-scale Lipschitz} if there exists a constant $\lambda>0$ and a constant $C\geq 0$ such that for all $x_1,x_2\in X,$ $\dist_Y(\phi(x_1),\phi(x_2))\leq \lambda \dist_X(x_1,x_2)+C.$

The map $\phi$ is  \textit{bornologous} if for every $\epsilon>0$ there exists $\delta >0$ such that for all $x_1,x_2\in X$, if  $\dist_X(x_1,x_2)<\epsilon$ then $\dist_Y(\phi(x_1),\phi(x_2))\leq \delta.$ 

The map $\phi$ is  \textit{coarse} if it is bornologous and the inverse image, under $\phi,$ of each bounded subset of $Y$, is a bounded set of $X.$

The maps $\phi,\phi' \colon X\to Y$ are \textit{close} to each other if there exists $C>0$ such that for all $x\in X,$ $\dist_Y(\phi(x),\phi'(x))<C.$ 

Let $\mathcal C$ be a class of maps. We say that $X$ and $Y$ are {\it $\mathcal C$-equivalent} if there exist maps $\phi \colon X\to Y$ and $\psi \colon Y\to X$ in the class $\mathcal C$ such that $\phi\circ \psi$ and $\psi \circ \phi$ are close to the identity. 

Two metric spaces are  $X$ and $Y$ {\it quasi-isometric} if they are large-scale Lipschitz equivalent.

A \textit{uniform spanning tree} of $\Gamma$ is a tree $T$ with a bijective map $\phi \colon VT\to V\Gamma$ such that $\phi$ and $\phi^{-1}$ are Lipschitz continuous.
\end{definitions}

The following result follows immediately from the definitions.
\begin{lemma}\label{L:mapclass}
$\phi$ Lipschitz continuous $\Rightarrow$ $\phi$ large-scale Lipschitz $\Rightarrow$  $\phi$ bornologous.\qed
\end{lemma}

\begin{remark}\label{R:finind}
Recall that the Cayley graphs of $G$ with respect to different generating sets are quasi-isometric. Hence the quasi-isometric equivalence class of $G$ is well defined. Recall also that a group $G$ is quasi-isometric to any finite index subgroup. Hence Cayley graphs of virtually free groups are quasi-isometric to trees. See \cite[Proposition 3.19]{GhysHarpe}. 

\end{remark}

\begin{definitions}
Let  $m,n\in [0\uparrow \infty[.$ 

An \textit{$m$-sequence} of length $n$ in $V\Gamma$ is a finite sequence $v_{[0\uparrow n]}$ in $V\Gamma$ such that $v_n=v_0$ and $\dist_\Gamma(v_i,v_{i+1})\leq m$ for $i\in[0\uparrow n-1].$

An $m$-sequence $v_{[0\uparrow n]}$ in $V\Gamma$ is \textit{$m$-reducible} if there exists $i,\in[1\uparrow n-1],$ such that $\dist_\Gamma(v_{i-1},v_{i+1})\leq m,$ in this event we can perform an \textit{$m$-reduction} of $v_{[1\uparrow n]}$ and obtain a new $m$-sequence $v_{[0\uparrow i-1]}\vee v_{[i+1\uparrow n]}.$

An $m$-sequence is \textit{$m$-triangulable} if it has length $\leq 3$ or if there exists an $m$-reduction that produces a sequence which is $m$-triangulable.

\end{definitions}

\begin{example}\label{Ex:0red}
Recall that a closed path $(v_{[0\uparrow n]},e_{[1\uparrow n]})$ is {\it reducible } if $e_1=e_n^{-1}$ or $e,e^{-1}$ is a subsequence of $e_{[1\uparrow n]}.$

Hence if a closed path $(v_{[0\uparrow n]},e_{[1\uparrow n]})$ of length $\geq 3$ is reducible then the $1$-sequence $v_{[0\uparrow n]}$ is $0$-reducible.
\end{example}

\begin{remark}
The term triangulation can be justified as follows. Let $v_{[0\uparrow n]}$ be a 1-sequence which is $m$-triangulable. For each $k\in [0\uparrow n],$ identify $v_k$ with the complex number $e^{2ki\pi/n},$ and the $1$-sequence with the regular polygon $P$ of vertices $v_{[0 \uparrow n]}.$ 
If $n\leq 3,$ then $P$ is a (maybe very degenerate) triangle.

Let  $v_{[0\uparrow p-1]}\vee v_{[p+1\uparrow n]}$ be the $m$-sequence obtained by an  $m$-reduction of the $m$-triangulation. Add a diagonal from $e^{2(p-1)\pi i/n}$ to $e^{2(p+1)\pi i/n}.$ This diagonal divides $P$ into a polygon with $n-1$ sides and a triangle. Repeating this process, the $m$-triangulation of $v_{[0\uparrow n]}$ gives a triangulation of $P$ by diagonals.
\begin{figure}[ht]
\begin{center}
\scalebox{0.6}{\includegraphics{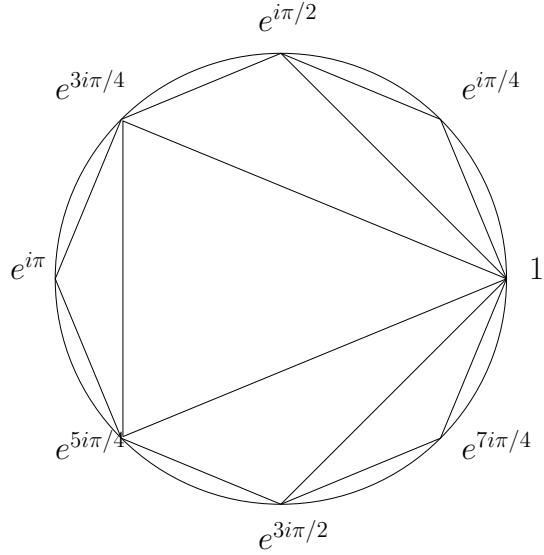}}
\caption{A triangulation in the complex plane.}
\label{fig:triangulation}
\end{center}
\end{figure}
\end{remark}

We still need some more definitions, however we state now the main theorem of this section. The proof is broken into a series of Lemmas, and definitions will be introduced as needed. 

\begin{theorem}\label{T:Eqcond}
Let $\Gamma$ be graph.
Then the following are equivalent:
\begin{enumerate}
\item[{\normalfont (B1).}] $\Gamma$ admits a uniform spanning tree.
\item[{\normalfont(B2).}] $\Gamma$ is quasi-isometric to a tree.
\item[{\normalfont(B3).}] $\Gamma$ is bornologous equivalent to a tree.
\item[{\normalfont(B4).}] There exists $m>1$ such that every $1$-sequence in $V\Gamma$ is $m$-triangulable.
\item[{\normalfont(B5).}] There exists $m>1$ such that for every $u,v,z\in V\Gamma$ and every path $\gamma$ from $u$ to $v$ the inequality \begin{equation*}\dist(z,\gamma)\leq \dfrac{1}{2}(\dist(u,z)+\dist(z,v)-\dist(u,v)+3m)\end{equation*} holds.
\item[{\normalfont(B6).}] There exists $m>1$ such that for every $z\in V\Gamma$ and every $n\in[0\uparrow \infty[,$ the boundaries of $\Gamma-\{g\in G : |g|\leq n\}$ have diameter less than or equal to $m.$
\item[{\normalfont(B7).}] $\Gamma$ admits a strong tree decomposition of uniformly bounded diameter.
\end{enumerate}
Moreover if $\Gamma$ is the Cayley graph of $G$ with respect to a finite generating set,
then (B1)-(B7) and the following are equivalent:
\begin{enumerate}
\item[{\normalfont (B8).}] The ends of $\Gamma$ have uniformly bounded diameter.
\item[{\normalfont (B9).}] There exists $m'>0$ such that for every $u,v\in V\Gamma,$ every path $\gamma$ from $u$ to $v$ and every $z$ in the geodesic from $u$ to $v$ the inequality $\dist(z,\gamma)\leq m'$ holds.
\item[{\normalfont (A3).}] $G$ is virtually free.
\item[{\normalfont (A4).}] There exists a finite set $S$ and an almost-right-$G$-invariant function $\alpha\in(G,S)$ such that for all $\beta=_a \alpha$ the right-stabilizer of $\beta$ is finite. 
\end{enumerate}
\end{theorem}
\begin{proof} 

Notice that the value of the constant $m$  may be different from one condition to another.

By Theorem~\ref{T:A} (A4)$\Rightarrow$(A3), by Remark~\ref{R:finind} (A3)$\Rightarrow$(B2). 
By Lemma~\ref{L:mapclass}, (B1) $\Rightarrow$ (B2) $\Rightarrow$ (B3).

The scheme of the rest of the proof is the following. First we close a cycle of implications from  (B1) to (B7) and back to (B1).

Then we show (B7)$\Rightarrow$ (B8)$\Rightarrow$ (B9)$\Rightarrow$ (A4). Notice that (B9) is a particular case of (B5), with $m'=3/2 m$. The reader interested in showing  the Muller and Schupp Theorem or that virtually free groups are quasi-isometric to trees can go from Lemma~\ref{L:triang->pathineq} directly to Lemma~\ref{L:pathineq2->A}.

We start now the proof of the Theorem. 
The following useful Lemma is a classic result of Nielsen \cite{Nielsen21}. See also \cite[I Proposition 2.2]{LyndonSchupp}.
\begin{lemma}\label{L:triangintree}
Let $T$ be tree. For all $m\in[0\uparrow \infty[,$ every $m$-sequence in $VT$ is $m$-triangulable.
\end{lemma}
\begin{proof}

We will argue by induction on the length of the $m$-sequence. By definition all $m$-sequences of length $\leq 3$ are $m$-triangulable. Let $v_{[0\uparrow n]}$ be an $m$-sequence in $VT$ with $n\geq 4$ and assume that any shorter $m$-sequence is $m$-triangulable.

Let $T'$ be the smallest connected subgraph of $T$ containing $v_{[0\uparrow n]}.$ Hence $T'$ is a locally finite tree and it can be identified with a subgraph of the Cayley graph of a free group $F$ with respect to a free generating set $X.$ We may think of $v_{[0\uparrow n]}$ as a sequence in $F.$ 

By the induction hypothesis, it is enough to show that $v_{[0\uparrow n]}$ is $m$-reducible

For $i\in[1\uparrow n],$ let $w_i=v_{i-1}v_{i}^{-1}.$ Hence $\Pi w_{[1\uparrow n]}=1.$ Notice that for $i\in[1\uparrow n],$ \begin{equation*}|w_i|=\dist_X(1,w_i)=\dist_T(v_{i-1},v_i)\leq m.\end{equation*}

For $i\in[2\uparrow n],$ let $a_i$ be the maximal common prefix of $w_i$ and $w_{i-1}^{-1}.$ Let $c_{i-1}=a_{i}^{-1}.$  Similarly, let $a_1=c_n^{-1}$ be the common prefix of $w_1$ and $w_{n}^{-1}.$ 
If for some $i\in[1\uparrow n],$ $|a_i|+|c_i| > |w_i|$, then either $|a_i|$ or $|c_i|$ is greater than half $|w_i|$. Suppose that $|a_i|>|w_i|/2.$ After cyclically permuting $v_{[0\uparrow n]},$ if necessary, we can assume that  $i\geq 2.$ Then $d_T(v_{i-2},v_i)=d_X(1,w_{i-1}w_i)\leq d_X(1,w_{i-1})\leq m$ and   $v_{[0\uparrow n]}$ is $m$-reducible.

We restrict then to the case $|a_i|+|c_i| \leq |w_i|,$ for all $i\in [0\uparrow n].$ Let $b_i=a_i^{-1}w_ic_i^{-1}.$

Hence \begin{equation*}1=\Pi w_{[1\uparrow n]}={a_1}^{-1}(\Pi b_{[1\uparrow n]})c_n^{-1}.\end{equation*}
By construction, there are no reductions on the product $\Pi b_{[1\uparrow n]},$ $b_i=1$ for $i\in[1\uparrow n].$

Let $j\in[2\uparrow n-2].$ 

If $|w_jw_{j+1}|\leq m$ then $\dist_T(v_{j-1},v_{j+1})<m$ and hence $v_{[0\uparrow n]}$ is $m$-reducible.

If $|w_jw_{j+1}|>m,$ as $|w_jw_{j+1}|=|a_j|+|c_{j+1}|>m$  either $|a_j|>m/2$ or $|c_{j+1}|>m/2.$ Suppose the former, the discussion for the latter is analogous. As $|a_j|>m/2,$ $|c_j|<m/2$ also as $|c_{j-1}|=|a_j|,$  $|a_{j-1}|<m/2.$ Hence  
\begin{equation*}m>|a_{j-1}|+|c_j|=\dist_X(1,w_{j-1}w_{j})=\dist_T(v_{j-2},v_{j}).\end{equation*}
and $v_{[0\uparrow n]}$ is $m$-reducible.

In the case $|c_{j+1}|>m/2$ we conclude that $\dist(v_{j},v_{j+2})\leq m$ and hence $v_{[0\uparrow n]}$ is $m$-reducible. 
\end{proof}

\begin{corollary}[(B3) $\Rightarrow$(B4)]\label{C:born->trian}
Let $T$ be a tree and suppose that $\Gamma$ is bornologous equivalent  to $T$. Then for every $k\in[1\uparrow \infty[,$ there exists $k'\geq 0$ such that every $k$-sequence in $V\Gamma$ is $k'$-triangulable. \hfill\qed
\end{corollary}
\begin{proof}
Let $\phi \colon V\Gamma\to VT$ and $\psi \colon VT\to V\Gamma$ be bornologous maps, such that $\psi\circ \phi$ is close to the identity. Let $R\in [1\uparrow \infty[$ such that for every $v\in V\Gamma,$ $\dist(v,\psi\circ \phi(v))<R$.

Let $k\in[0\uparrow \infty[.$  There exists $m\in[0\uparrow \infty[$ such that for all $u,v\in V\Gamma,$ if $\dist_{\Gamma}(u,v)<k$ then \begin{equation*}\dist_T(\phi(u),\phi(v))<m.\end{equation*}
There exists $k_1>0$ such that if $\dist_T(\phi(u),\phi(v))<m$ then $\dist_\Gamma(\psi\circ\phi(u),\psi\circ\phi(v))<k_1,$ and hence $\dist_\Gamma(u,v)\leq 2R+k_1.$

Let $v_{[0\uparrow n]}$ be a $k$-sequence in $V\Gamma$. Hence $\phi(v_{[0\uparrow n]})$ is an $m$-sequence in $VT$ and, by Lemma~\ref{L:triangintree}, it is $m$-triangulable. As each $m$-reduction of $\phi(v_{[0\uparrow n]})$ gives a $(2R+k_1)$-reduction of $v_{[0\uparrow n]}$, it follows that $v_{[0\uparrow n]}$ is $(2R+k_1)$-triangulable.
\end{proof}

\begin{lemma}[(B4) $\Rightarrow$(B5)]\label{L:triang->pathineq}
Let $m>1.$ Let $\Gamma$ be a graph such that every 1-sequence in $V\Gamma$ is $m$-triangulable. Then for every $u,v,z\in V\Gamma$ and any path $x_{[0\uparrow n]}$ in $V\Gamma$ from $u$ to $v,$ the following inequality holds:
\begin{equation}\label{eq:pathineq}
\dist_\Gamma(z,x_{[0\uparrow n]})\leq \dfrac{1}{2}(\dist_\Gamma(z,u)+\dist_\Gamma(z,v)-\dist_\Gamma(u,v)+3m).
\end{equation}
\end{lemma}
\begin{proof}
Let $u_{[0\uparrow \dist(z,u)]}$ and $v_{[\dist(z,v)\downarrow 0]}$ be  geodesic paths in $V\Gamma$ from $z$ to $u$ and from $v$ to $z,$ respectively.
Then $u_{[0\uparrow \dist(z,u)]}\vee x_{[1\uparrow n-1]}\vee v_{[\dist(z,v)\downarrow 0]}$ is a $1$-sequence and hence it is $m$-triangulable. Notice that $u_0=v_0=z.$

If $\dist(u,v)\leq 1,$ as $m>1,$ the inequality~\eqref{eq:pathineq} holds. Thus, we may assume that $\dist(u,v)>1.$

There exists a sequence $s_{[0\uparrow N]}$ of $m$-sequences such that $s_0=u_{[0\uparrow \dist(z,u)]}\vee x_{[1\uparrow n-1]}\vee v_{[\dist(z,v)\downarrow 0]},$  for $i\in[1\uparrow N]$ each $s_i$ is an $m$-reduction of $s_{i-1},$ and $s_N$ has length $3.$

As $u_0=v_0=z$ are the first and last terms in all $s_i,$ $i\in[0\uparrow N],$ there exists $j\in[0\uparrow N]$ such that  $s_j$ only contains one term of $x_{[1\uparrow n-1]},$ say $x_b.$ 

If $s_j$ has the form $(u_0,x_b,\dots,v_0)$ or $(u_0,\dots, x_b,v_0)$ then $\dist(z,x_{[0\uparrow n]})\leq m$ and the inequality~\eqref{eq:pathineq} hold. Thus, it remains the case when  $s_j$ has the form \begin{equation*}(u_0,\dots,u_a,x_b,v_c,\dots,v_0)\end{equation*} where  $a\in]0\uparrow \dist(z,u)],$ $c\in[\dist(z,v)\downarrow 0[,$ and $s_{j+1}$ is obtained from $s_j$ by suppressing $x_b.$

\begin{figure}[ht]\label{figabc}

\begin{center}
\scalebox{0.5}{\includegraphics{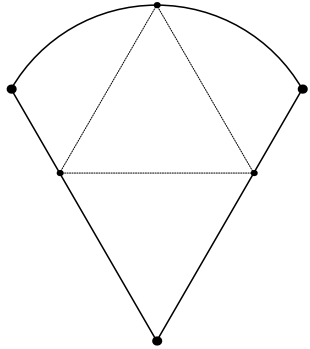}}
\end{center}
\caption{Triangle $u_a,x_b,v_c$ and the $1$-sequence $s_0.$}

\end{figure}
In this event $(u_a,x_b,v_c,u_a)$ is an $m$-sequence.  Then
\begin{eqnarray*}
\dist(u,v)&\leq&\dist(u,u_a)+m+\dist(v_c,v).\\
\dist(z,x_{[0\uparrow n]})&\leq&\dist(z,u_a)+m.\\
\dist(z,x_{[0\uparrow n]})&\leq&\dist(z,v_c)+m.\\
\end{eqnarray*}
Summing up the three inequalities we get
\begin{equation*}\dist(u,v)+2\dist(z,x_{[0\uparrow n]})\leq \dist(z,u) + \dist(z,v) + 3m\end{equation*}
which is equivalent to \eqref{eq:pathineq}.
\end{proof}
\begin{remark}
Recall that the {\it Gromov product } of $x,y\in V\Gamma$ based at $z$ is
\begin{equation*}(x\mid y)_z=\dfrac{1}{2}(\dist(z,x)+\dist(z,y)-\dist(x,y)).\end{equation*}
Hence \eqref{eq:pathineq} can be reformulated as for all $x,y,z\in V\Gamma$ and all paths $\gamma$ from $x$ to $y,$ $\dist(z,\gamma)\leq (x\mid y)_z+\frac{3}{2}m.$
\end{remark}

As we mentioned, readers not interested in (B6)-(B8) can go directly to Lemma~\ref{L:pathineq2->A}.

\begin{definitions}
Let $S\subset \Gamma$ and $z\in V\Gamma.$ The {\it boundary} of $S$ is \begin{equation*}\partial S=\{v\in VS : \exists u\in V\Gamma-VS\text{ such that }\dist(u,v)=1\}.\end{equation*}
The {\it ball of center $z$ and radius} $n,$ written $\ball_z(n),$ is the complete subgraph of $\Gamma,$ having  all the vertices $v\in V\Gamma$ at distance at most $n$ from $z$ as vertex set.
\end{definitions}

\begin{lemma}[(B5)$\Rightarrow$(B6)]\label{L:pathineq->boundaries}
Let  $m>1.$ Suppose that for every $u,v,z\in V\Gamma$ and any path  $x_{[0\uparrow n]}$ in $V\Gamma$ from $u$ to $v$ the inequality~\eqref{eq:pathineq} hold. Then for all $z\in V\Gamma$ and all $n\in [0\uparrow \infty[,$ all connected components $C$ of $\Gamma-\ball_z(n)$ satisfy  $\diam(\partial C)\leq 3m/2.$
\end{lemma}
\begin{proof}
Fix $z\in V\Gamma$ and $n\in [0\uparrow \infty[$ and let $u,v\in V\Gamma$ be in the boundary of a connected component of $\Gamma-\ball_z(n).$ Let $x_{[0\uparrow l]}$ be a path in $V\Gamma$ from $u$ to $v$ contained in $\Gamma-\ball_z(n).$ Then by \eqref{eq:pathineq}, $n=\dist(z,x_{[0\uparrow l]})\leq 1/2(2n-\dist(u,v)+3m)$ and hence $\dist(u,v)\leq 3m/2.$
\end{proof}
\begin{definition} \label{def:rgraphpartition}
Let $(X,\dist)$ be a metric space,  $\mathcal P$  a partition of $X$ and $r>0$.
The {\it $r$-graph associated to $\mathcal P$} is an unoriented graph $\Delta$ with vertex set $\mathcal P$ and an edge between $v,u\in \mathcal P$ if and only if $\dist(u,v)<r$.

A {\it strong tree decomposition} of $V\Gamma$ is a partition $\mathcal P$ of $V\Gamma$ such that 
the $1$-graph associated to $\mathcal P$ is a tree.

A strong tree decomposition $\mathcal P$ has {\it uniformly bounded diameter} (resp. {\it finite width}) if there exists a constant $K>0$ such that for all $S\in \mathcal P,$ $\diam_\Gamma(S)<K$  (resp. $|S|<K$). 
\end{definition}
\begin{lemma}\label{L:finwith}
If there is a strong tree decomposition $\mathcal P$ of $\Gamma$ with uniformly bounded diameter and the degrees of the vertices of $\Gamma$ are uniformly bounded, the strong tree decomposition also has finite width.
\end{lemma}
\begin{proof}
Let $k>0$ such that for all $S\in \mathcal P,$ $\diam_\Gamma(S)<K.$
Let $k_1$ be a uniform bound of the degrees of the vertices of $\Gamma.$
For $S\in \mathcal P,$ let $Y$ be a minimal subtree of $\Gamma$ containing all vertices of $S.$ A ball of radius $k$ in a regular tree of degree $k_1$ has at most $k_1^{k}$ vertices, it follows that $Y$ has at most $k_1^{k}$ vertices and hence, $|S|<k_1^{k}.$
\end{proof}

\begin{lemma}[(B6)$\Rightarrow$(B7)]\label{L:treedecom}
Let $\Gamma$ be a connected graph and $z\in V\Gamma.$  Suppose that for all $n\in [0\uparrow \infty[,$ the boundaries of the connected components of $\Gamma-\ball_z(n)$ have diameter  $\leq k.$ Then there exists a strong tree decomposition of $\Gamma$ with uniformly bounded diameter.

Moreover, if the degrees of the vertices of $\Gamma$ are uniformly bounded by $k_1>0,$ the strong tree decomposition also has  finite width.
\end{lemma}

\begin{proof}

Let \begin{equation*}\mathcal{P}=\{z\}\cup \{\text{boundaries of connected components of } \Gamma-\ball_z(n): n\in[1\uparrow \infty[\}.\end{equation*} 

Hence, for all $S\in \mathcal P,$ $\diam_\Gamma (S)<k.$

Let $\Delta$ be the $1$-graph associated to the partition $\mathcal P.$ To show that $\Delta$ is a tree, it is enough to show that a closed path $(S_{[0\uparrow l]},e_{[1\uparrow l]})$  in $\Delta$ is reducible.

Without loss of generality, assume that $l>1.$ By construction, there are no edges joining a vertex with itself, hence, for $i\in[1\uparrow l],$ $S_{i-1}\neq S_{i}$. There exists $m\in[1\uparrow l]$ such that for all $i\in [0\uparrow l],$ $\dist(z, S_m)\geq \dist(z,S_i).$  For some $t>0,$ $S_{m-1}$ and $S_{m+1}$ are the boundaries  of the connected component of $\Gamma-\ball_z(t)$ containing $S_m.$ Hence $S_{m-1}=S_{m+1}$ and  the path is reducible. 

\begin{figure}[ht]\label{figstree}
\begin{center}
\scalebox{0.5}{\includegraphics{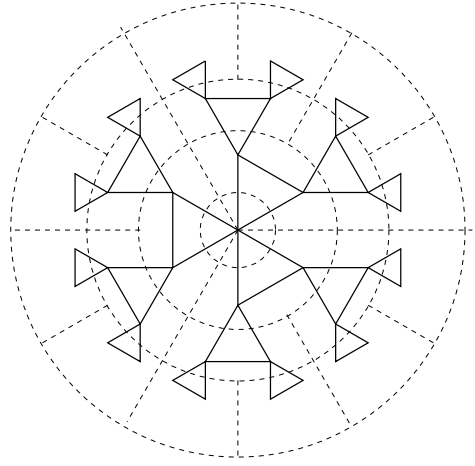}}
\end{center}
\caption{A tree decomposition.}
\end{figure}
If the degrees of the vertices of $\Gamma$ are unifomly bounded, then, by Lemma~\ref{L:finwith}, the strong tree decomposition also has finite width.
\end{proof}

\begin{lemma}[(B7)$\Rightarrow$(B1)]\label{L:std->usp}
Suppose that $\Gamma$ admits a strong tree decomposition of uniformly bounded diameter. Then $\Gamma$ admits a uniform spanning tree. 
\end{lemma}
\begin{proof}
Let $\mathcal P$ be a partition of $V\Gamma$ into sets of diameter $\leq k$  such that the $1$-graph associated to $\mathcal P$ is a tree. For each $S\in \mathcal P,$ take $v_S\in V\Gamma$ with $v_S\in S.$ 

Let $T$ be the tree with vertex set $V\Gamma,$ such that for all $S\in \mathcal P$ and all  $v\in S-\{v_S\}$ there is an edge joining $v$ and $v_S,$ and for all $S,R\in \mathcal{P}$ with $\dist(R,S)=1$ there is an edge joining $v_S$ and $v_R.$ 

Let $u,v\in V\Gamma$, $u\neq v,$ $u\in S,$ $v\in R,$ and $S,R\in \mathcal P.$

If $R=S$ then $\dist_T(u,v)\leq 2\leq 2 \dist_\Gamma(u,v)\leq 2k\dist_T(u,v).$ 

If $R\neq S$ then 
\begin{align*}
\dist_T(u,v)&\leq 2+\dist_T(v_S,v_R)\\%
&\leq 2+\dist_\Gamma(R,S)\\%
&\leq 2+\dist_\Gamma(u,v)\\%
&\leq 2\dist_\Gamma(u,v)+\dist_\Gamma(u,v)\\%
&\leq 3\dist_\Gamma(u,v)\end{align*}

Also $\dist_\Gamma(u,v)\leq (2k+1)\dist_T(u,v).$
Hence for all $u,v\in V\Gamma,$ \begin{equation*}\dfrac{1}{3}\dist_T(u,v)\leq \dist_\Gamma(u,v)\leq (2k+1)\dist_T(u,v).\end{equation*} 
\end{proof}

Notice that we have proved the equivalence between (B1)-(B7).

\begin{remark}\label{R:treedecom}
There exists the more general concept of a {\it tree decomposition} of a graph, which we will not use here. This concept is important in the theory of graph minors. In \cite{KuskeLohrey}, Kuske and Lohrey characterize virtual freeness in terms of the tree decompositions.

By Theorem 3.4 and Theorem 3.7 of \cite{KuskeLohrey}, a vertex-transitive locally-finite graph admits a strong tree decomposition of finite width and uniformly bounded diameter if and only if $\Gamma$ has a tree decomposition of finite width.

Hence, for Cayley graphs, all these concepts turn out to be equivalent.
\end{remark}

\begin{definitions}
Let $\Gamma$ be a locally finite graph. Let $\mathfrak B \Gamma=\{B\subset V\Gamma: |\partial B|<\infty\}.$  Let $\mathfrak P\Gamma$ be the set of infinite paths in $V\Gamma$ which do not repeat a vertex.

For all $B\in \mathfrak B \Gamma$ and $p\in \mathfrak P \Gamma,$ as $|\partial B|<\infty,$ either $p\subset_a B$ or $p\subset_a \Gamma-B.$ For each $p\in \mathfrak P \Gamma,$ we define a map $[p]\colon \mathfrak B \Gamma \to \{0,1\}$ as follows:  for $B\in \mathfrak B \Gamma$ we set $[p](B)=1$ if $p\subset_a B$ and $[p](B)=0$ otherwise. There is a natural equivalence relation $\sim$ in  $\mathfrak P \Gamma$ given by $p\sim p'$ if and only if $[p]=[p'].$

The space $\mathfrak E \Gamma\coloneqq \mathfrak P\Gamma\slash\hspace{-1mm}\sim \,= \{[p] : p \in \mathfrak P \Gamma  \}$ is called {\it the space of ends of $\Gamma$}. For $B\in \mathfrak B \Gamma,$ let \begin{equation*}\ol B=B\cup\{[p]\in \mathfrak E\Gamma: [p](B)=1\}\subseteq \mathfrak B \Gamma \cup \mathfrak E \Gamma\end{equation*} and let $\ol{\mathfrak B \Gamma}=\{\ol B: B\in \mathfrak B \Gamma\}.$ Then $V\Gamma\cup \mathfrak E \Gamma$ with the topology defined by $\ol{\mathfrak B\Gamma}$ is a totally disconnected compact Hausdorff topological space. See \cite[IV.6]{DicksDunwoody89}, for details.

Let $k>0.$ An {\it end $[p]$ of $\mathfrak E \Gamma$ has diameter $<k$} if there exists a sequence $B_{[0\uparrow \infty[}$ in $\mathfrak B \Gamma$ such that for all $n\in[0\uparrow \infty[,$ $|\partial B_n|<k,$ $B_n\subset B_{n+1}$ and  $\cap \ol{B}_{[0\uparrow \infty[}=[p].$ The sequence $B_{[0\uparrow \infty[}$ is {\it a sequence of diameter $<k$ converging to $[p]$}.

The ends of $\Gamma$ have {\it uniformly bounded diameter} if there exists $k>0$ such that every $[p]\in \mathfrak E \Gamma$ has diameter $<k.$
\end{definitions}

\begin{lemma}[(B7)$\Rightarrow$(B8)]\label{L:std->ends}
If the vertices of $\Gamma$ have uniformly bounded degree and $\Gamma$ has a strong tree decomposition of uniformly bounded \hyphenation{di-a-me-ter} diameter, then the ends of $\Gamma$ have uniformly bounded diameter.
\end{lemma}
\begin{proof}
By Lemma~\ref{L:treedecom}, there exist $k$ and $k_1$ and a strong tree decomposition $\mathcal P$ of $\Gamma$ such that for all $S\in \mathcal P,$ $\diam S<k$ and $|S|<k_1.$

We may consider $\mathcal P$ as the vertex set of a locally finite tree $T.$ 

Let $p=v_{[0\uparrow \infty]}\in  \mathfrak P \Gamma.$  The path $p$ can be viewed as an infinite sequence $v_{[0\uparrow \infty[}$ in $VT=\mathcal P$ with $\dist(v_i,v_{i+1})\leq 1$ which repeats a vertex at most $k_1$ times. Hence for $B\in \mathfrak B T,$ either $p\subset_a B$ or $p\subset_a VT-B.$ 

We construct a sequence $B_{[0\uparrow  \infty]}$ in $\mathfrak B\Gamma$ such that $B_{i+1}\subset B_i$ and $[p](B_i)=1$ for $i\in[0\uparrow \infty].$ Let $B_0$ be the connected component of $T-\{v_0\}$ such that $[p](B_0)=1.$ To construct $B_{i+1}$ from $B_{i},$ we consider the first $j\in[0\uparrow \infty[$ such that $v_{k}\in B_i$ for $k\in[j\uparrow \infty[$. As $[p](B_i)=1$ such a $j$ exists. Then there exists a connected component $C$ of $T-\{v_j\}$ strictly contained in $B_i$ such that $[p](C)=1$. Let $B_{i+1}=C.$

Hence for $i\in [0\uparrow \infty[,$ $|\partial B_i|<k,$ $B_i\subset B_{i+1}$ and $[p]\in \cap \ol{B}_{[1\uparrow \infty[}.$ 

Let $B\in \mathfrak B\Gamma.$  The set $B$ is almost equal to a subforest of $T$ with finitely many components. Suppose that $(\cap \ol{B}_{[0\uparrow \infty[})\cap \ol{B}\neq \emptyset,$ hence for a large enough $n,$ $B_n$ is contained in a subforest of $B.$

Let $q \in \mathfrak P \Gamma$ such that $[q]\in \cap \ol{B}_{[0\uparrow \infty[}.$ If $[q](B)=1,$ then $(\cap \ol{B}_{[0\uparrow \infty[})\cap \ol{B}\neq \emptyset$ and $B$ contains some $B_n,$ therefore $[p](B)=1.$ Similarly, if $[q](B)=0$ then $[p](B)=0.$ Hence $[p](B)=[q](B)$ for all $B\in \mathfrak B \Gamma,$ thus $p\sim q$ and $[p]=\cap \ol{B}_{[0\uparrow \infty[}.$
\end{proof}

\begin{lemma}[(B8)$\Rightarrow$(B9)]\label{L:ends->pathineq2}
Let $\Gamma$ be a locally finite vertex-transitive graph. Suppose that the ends of $\Gamma$ have uniformly bounded diameter. Then there exists $m>0$ such that for all $u,v\in V\Gamma,$ any $\eta$  geodesic from $u$ to $v,$ any $z\in\eta$ and any path $\gamma$ from $u$ to $v,$ $\dist(z,\gamma)\leq m.$
\end{lemma}

\begin{proof}
Let $k>0$ be a uniform bound of the diameter on the ends of $\Gamma.$ Let $t\in V\Gamma.$ For $[p]\in \mathfrak E\Gamma,$ let $B^{[p]}_{[0\uparrow \infty[}$ be a sequence of diameter $<k$ converging to $[p].$ The set 
\begin{equation*}\mathcal C=\{ B_n^{[p]} : [p]\in\mathfrak E \Gamma, n\in [0\uparrow \infty[, \dist(t,B_n^{[p]})>k \}\end{equation*} 
is an open covering of $\mathfrak E \Gamma.$ As $\mathfrak E\Gamma$ is compact, there exist $B_0,\dots, B_n\in \mathcal C$ such that $\mathfrak E\Gamma\subset \cup \ol{B}_{[0\uparrow n]}.$

Let $S=V\Gamma- \cup B_{[0\uparrow n]}.$ Then $S$ is finite and  $t\in S.$ Let $m>0$ such that
\begin{equation}\label{eq:diamS}
\dist(\partial B_i, t)<m \text{ for all }i\in[0\uparrow n] 
\end{equation}
 and $\diam(S)<m.$ 

Let $u,v,z,\eta$ and $\gamma,$ be as in the hypothesis. As $\Gamma$ is vertex transitive, we may assume that $z=t.$

If $u\in S$ or $v\in S,$ then  $\dist(z,\gamma)\leq m.$

If $u\notin  S$ and $v\notin S,$ then $u\in B_i$ and $v\in B_j$ for some $i,j\in [0\uparrow n].$ If $i=j,$ then there exist $u',v'\in \partial B_i$ such that \begin{equation*}\dist(u,z)=\dist(u,u')+\dist(u',z)\end{equation*} and \begin{equation*}\dist(v,z)=\dist(v,v')+\dist(v',z).\end{equation*} 
As $\diam(\partial B_i)<k,$
\begin{equation}\label{eq:contra1}
\dist(u,u')+k+\dist(v',v)\geq \dist(u,u')+\dist(u',v')+\dist(v',v)\geq \dist(u,v).
\end{equation}
As $z$ is in the geodesic $\eta,$ $\dist(u,v)=\dist(u,z)+\dist(z,v).$
By construction of $S,$ $\dist(u',z)>k$ and $\dist(v',z)>k,$ hence
\begin{align}
\label{eq:contra2}
\dist(u,v)&=\dist(u,z)+\dist(z,v)\nonumber \\%
&=\dist(u,u')+\dist(u',z)+\dist(z,v')+\dist(v',v) \nonumber\\%
 &\geq \dist(u,u')+2k+\dist(v,v').
\end{align}
Combining \eqref{eq:contra1} and \eqref{eq:contra2}, $\dist(u,u')+k+\dist(v,v')\geq \dist(u,u')+2k+\dist(v,v'),$ but this is impossible since $k>0.$

Hence $i\neq j$ and then path $\gamma$ goes through $\partial B_i.$ Hence, by \eqref{eq:diamS}, $\dist(z,\gamma)<m.$
\end{proof}

Now we connect the graph conditions using Theorem~\ref{T:Dun}.
\begin{lemma}[(B9)$\Rightarrow$(A4)]\label{L:pathineq2->A}
Let $\Gamma$ be a Cayley graph of a group $G$ with respect to a finite generating set $X.$ Let $m>0$ such that for all $u,v\in V\Gamma,$ any  geodesic $\eta$ from $u$ to $v,$ any $z\in\eta$ and any path $\gamma$ from $u$ to $v,$  $\dist(z,\gamma)\leq m.$

Then there exists a finite set $S$ and an almost-right-$G$-invariant function $\alpha\in (G,S)$ such that for all $\beta=_a \alpha,$ the right-stabilizer of $\beta$ is finite. 
\end{lemma}
\begin{proof}
If $G$ is finite, the result is trivial. So, we assume that $G$ is infinite.
Let $S$ be the partition of $G$ formed by the (vertex set of the) connected components of $\Gamma-\ball_1(m)$ and (the vertex set of) $\ball_1(m).$ Let $\alpha\colon G\to S$ be the function satisfying  $g\in \alpha(g)$ for all $g\in G.$

Fix $g\in G$, and let $x_{[1\uparrow|g|]}$ be a sequence in $X^{\pm 1}$ such that $\Pi x_{[1\uparrow |g|]}=g$. Let  $h\in G$ such that $|h|>|g|+m.$ Then $|h\Pi x_{[1\uparrow i]}|>m$ for all $i\in[1\uparrow |g|]$ and hence $h$ and $hg$ lie in the same connected component of $\Gamma-\ball_1(m).$ Thus  $\alpha(h)g=\alpha(h).$ This shows that $\alpha$ is almost-right-$G$-invariant.

Let $\beta=_a\alpha.$ There exists $n\in [0\uparrow \infty[$ such that $\alpha$ and $\beta$ agree on $\{g\in G : |g|>n\}.$ 

Let $g\in G$ such that $|g|>2n+2,$ and $g=\Pi x_{[1\uparrow |g|]}.$ Let $h=\Pi x_{[1\uparrow n+1]}.$  Hence, as the geodesic from $h^{-1}$ to $h^{-1}g$ goes through $1,$ any path from $h^{-1}$ to $h^{-1}g$ goes through $\ball_1(m)$ and hence $h^{-1}$ and $h^{-1}g$ lie in different connected components of $\Gamma-\ball_1(m).$ Therefore \begin{equation*}\beta(h^{-1}g)=\alpha(h^{-1}g)\neq \alpha(h^{-1})=\beta(h^{-1})\end{equation*} and hence $g$ does not stabilize $\beta.$ Thus $G_\beta$ is finite.
\end{proof}
Compare the proof of the  Lemma with the proof of (A3) implies (A4) in Theorem~\ref{T:A}.

This completes the proof of Theorem~\ref{T:Eqcond}
\end{proof}

\begin{Rems}
The equivalence between (B5) and (A3) is due to Muller and Schupp \cite{MullerSchupp83}. By \cite[Theorem I]{MullerSchupp83}, a group is context-free if and only if (B5) holds.
The equivalence between (B2) and (A3) is attributed to Gromov and appears explicitly in \cite[Thm 7.19]{GhysHarpe}.

The equivalence between (B1), (B4) and (B8) appears in the  article of Woess \cite{Woess}. Woess uses the result of Muller and Schupp to conclude the equivalence of these conditions with (A3).

The equivalence between  (B7) and (A3) is due to Kuske and Lohrey \cite{KuskeLohrey}. See Remark~\ref{R:treedecom}.
\end{Rems}

\begin{Rev}\label{Rev:lang}
The original result of Muller and Schupp \cite{MullerSchupp83} is about context-free groups. For completeness we recall the definition.
 
By $X^*$ we denote the free monoid freely generated by a set $X.$ A subset $L$ of $X^*$ is {\it a language over X}. We identify $X^*$ with the set of sequences in $X$ with the operation given by concatenation.

Usually the languages are classified either by the type of grammars that generate the language or by the type of automatons that recognize the sequences in the language. We follow the second approach.

 A language $L$ is {\it regular}, if it is accepted by a finite-state automaton.

 A {\it finite-state automaton} $A$ over $X$ is a quadruple $A=(V, v,\lambda, S)$ where $V$ is a finite set, called the set of states, $v\in V$ is the initial state, $\lambda\colon V\times X \to V$ is a function and $S$ is a subset of $V$, called the set of accepting states. It is useful to think $A$ as an oriented labelled graph with vertex set $V$ and edges labelled by $x\in X$ starting at $u$ and ending at $\lambda(u,x)$ for $u\in V.$
 
 A sequence $x_{[1\uparrow n]}$ in $X$ is recognized by an automaton if there exists a sequence $v_{[0\uparrow n]}$ in $V$  such that $v_0=v,$ $v_{i}=\lambda(v_{i-1},x_i)$ for $i\in[1\uparrow n]$ and $v_n\in S.$

A language $L$ is {\it context-free}, if it is accepted by a push-down automaton. It can also be defined as a language generated by a context-free grammar.

Roughly speaking, a push-down automaton is a finite-state automaton with a memory stack.  Now the transition function $\lambda,$  depends on the state and the top symbol of the stack, and outputs a new state and modify the top of the stack, by deleting, adding or keeping the top symbol.

Formally,  a {\it push-down automaton} $A$ over $X$ is a 6-tuple $A=(V,v,\lambda, S, \Sigma, \alpha),$ where $V,v$ and $S$ are as for the finite-state automaton, $\Sigma$ is a finite set and $\alpha \in \Sigma,$ and the function $\lambda\colon V\times X\times \Sigma^*\to V\times \Sigma^*$ satisfies that if $(u,\beta)=\lambda(w,x,\gamma)$ then $\beta$ is  either equal to $\gamma$ or obtained from $\gamma$ by either adding a symbol at the end of the sequence, or deleting the last symbol of the sequence.

A sequence $x_{[1\uparrow n]}$ in $X$ is recognized by the  automaton if there exist  sequences $v_{[0\uparrow n]}$ in $V$ and $\alpha_{[0\uparrow n]}$ in $\Sigma^*$  such that $(v_0,\alpha_0)=(v,\alpha),$  $(v_{i},\alpha_i)=\lambda(v_{i-1},x_i,\alpha_{i-1})$ for $i\in[1\uparrow n]$ and $v_n\in S.$

Given a group $G$ and a finite generating set $X$ of $G,$ we say that $(G,X)$ is regular (respectively context-free) if the set $L(G,X)$ of the sequences $x_{[1\uparrow n]}\in(X^{\pm 1})^*$ such that  $\Pi x_{[1\uparrow n]}=1$ in $G$ is a regular (respectively context-free) language. It can be shown that if $X'$ is another finite generating set of $G$ and  $L(G,X)$ is regular (respectively context-free), then $L(G,X')$ is also regular (respectively context-free).

The regular groups are exactly the finite ones. This is a classical result of  Anisimov.


The Muller-Shupp theorem states that context-free groups are exactly the virtually free ones. 

\end{Rev}

\begin{Rems}\label{Rem:otherchar}
There exist other characterizations of virtual-freeness that involve language conditions.

For example, in \cite{BridsonGilman}, Bridson and Gilman show that it is possible to put more structure on the spanning tree of (B1) of Theorem~\ref{T:Eqcond}. In particular they show that a group $G$ is virtually free if and only if there exists a finite generating set $X$ of $G$ such that $G$ admits a regular broomlike combing with respect to $X.$

We remark that the result holds for any finite generating set.

Let $\mu \colon X^*\to G$ be a surjective monoid morphism. A subset $L$ of $X^*$ is {\it a combing} with respect to $X$ if the restriction of $\mu$ to $L$ is a bijection with $G.$ The combing is {\it broomlike} if there exists a constant $k$ such that for all $u,v\in L$ where $\dist_{G}(\mu(u),\mu(v))=1,$ there exist $u',v',w\in X^*$ such that $u=wu',$ $v=wv'$ and $|u'|_{X}+|v'|_{X}\leq k.$ 

The subset $L$ is a {\it regular broomlike combing} with respect to $X$ for $G,$ if $L$ is a regular language and a broomlike combing with respect to $X$ for $G.$

The monoid $X^*$ has a natural tree structure, and hence $L$ inherits a tree structure. It is not hard to show that $L$ with the bijection $\mu\colon L \to G$ is a uniform spanning tree.

Another result in the same spirit, can be found in \cite{GilmanHermillerHoltRees}, where Gilman, Hermiller, Holt and Rees proved the following: a group $G$ is virtually free, if and only if there exists a finite generating set $X$ and a constant $m,$ such that every $m$-locally geodesic sequence in $G$ is geodesic.

We remark that this result does not hold in general for every finite generating set.

Let  $x_{[1\uparrow n]}$ be a sequence in $X^{\pm 1}.$  It is  is geodesic if $|\Pi x_{[1\uparrow n]}|=n.$ The sequence is a  $m$-locally geodesic sequence, if any subsequence of length $m$ is geodesic.

Again this characterization can be partially connected to the present work. Using the ideas of \cite[Proposition 1]{GilmanHermillerHoltRees}, it is not hard to show that if every $m$-locally geodesic sequence is geodesic, then every $1$-sequence is $(4m+1)$-triangulable.
\end{Rems}

\section{Asymptotic dimension}\label{S:adim}
In this section we give an elementary introduction to groups of asymptotic dimension 1. A general reference about the topic is \cite[Chapter 9]{Roe}. 

The concept of asymptotic dimension was introduced by Gromov in \cite{Gromov91}. He proposed two different definitions of what asymptotic dimension should be, which turn out to be equivalent. See \cite[Theorem 9.9]{Roe}. We will use the following.
\begin{definition}\label{D:adim}
A metric space $(X,\dist)$ has {\it asymptotic dimension} at most $n,$ if for every $m>0,$ there exists a partition $\cal P$ of $X$ into sets of uniformly bounded diameter such that the $m$-graph of $\cal P$ is $(n+1)$-colorable, where the $m$-graph of $\cal P$ has vertex set $\cal P$ and two elements of $\cal P$ are joined by an edge if they are at distance at most $m$. 

We write $\adim \Gamma\leq n$ if the asymptotic dimension of $\Gamma$ is at most $n,$ and  $\adim \Gamma=n$ if $\adim \Gamma\leq n$ and  $\adim \Gamma\nleq n-1.$
\end{definition}
The asymptotic dimension is invariant under coarse equivalence. Hence, the asymptotic dimension of a finitely generable group is just the asymptotic dimension of its Cayley graph with respect to a finite generating set.

The  following observation is straightforward from Definition~\ref{D:adim}. 
\begin{observation}\label{O:defadim}
The Cayley graph $\Gamma$ of a finitely generated group $G$ has asymptotic dimension at most $n,$ if for every $m>0$ there exist two functions $z,r\colon G\to G$ such that $g=z(g)r(g),$ $r(G)$ is finite and there is a function $\col\colon z(G)\to [0\uparrow n]$ such that for $z(g_1)\neq z(g_2)$ if $\col(z(g_1))=\col(z(g_2))$ then $\dist(g_1,g_2)>m.$ \qed
\end{observation}

\begin{remark}
If $\adim G=0$ then $G$ is finite.
\end{remark}

\begin{lemma}\label{L:adimtree}
Let $T$ be a tree with the degrees of the vertices uniformly bounded. Then, $\adim T\leq 1.$ 
\end{lemma}
\begin{proof}
Let $m>0$ and fix a vertex $z\in VT.$ 
Let \begin{equation*}A_{n,i}=\{v\in VT : nm\leq \dist(z,v)< (n+i)m\}.\end{equation*}
Let $\mathcal P$ be the partition of $VT$ consisting of connected components of $A_{n,2}$ intersected with $A_{n+1,1}$ for all $n\in [-1\uparrow \infty[.$ 

Since the degrees of the vertices are uniformly bounded, there exists $R>0$ such that $\diam S<R$ for all $S\in \mathcal P.$

For $S\in \mathcal P,$  let $\col(S)=0$ if $S\subset A_{2k,1}$ for some $k\in [0\uparrow \infty[,$ and $\col(S)=1$ otherwise.

\begin{figure}[ht]\label{fig:coronas}
\begin{center}
\scalebox{1}{\includegraphics{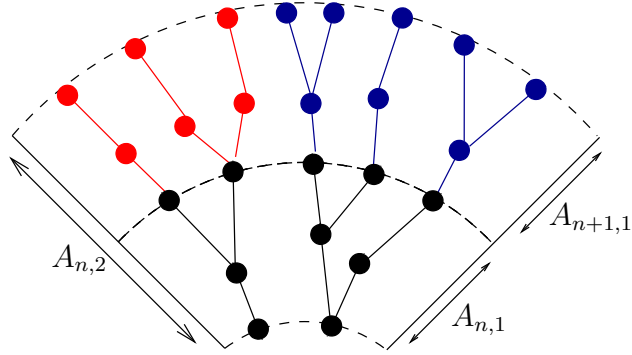}}
\end{center}
\caption{The coronas $A_{n,i}$ and the partition.}
\end{figure}

Let $S,S'\in \mathcal P$ with $\col(S)=\col(S'),$ $S\neq S'.$ Then $S\subset A_{2k+\varepsilon,1}$ and $S'\in A_{2k'+\varepsilon,1},$ with $\varepsilon\in \{0,1\}.$

If $k\neq k'$ then $\dist(S,S')\geq \dist(A_{2k+\varepsilon,1},A_{2k'+\varepsilon,1})\geq 2m.$
If $k=k'$ then $S$ and $S'$ lie in different connected components of $A_{2k+\varepsilon-1,2},$ and hence any path joining $S$ and $S'$ must go through $\{v\in VT:\dist(z,v)< (2k+\varepsilon-1)m\}$ and hence $\dist(S,S')\geq 2m.$
\end{proof}

\begin{corollary}
Let $G$ be a finitely generable group. If $G$ is virtually free, then $\adim G\leq 1.$
\end{corollary}
\begin{proof}
Let $H$ be a finite-index, free subgroup of $G.$ Since $H$ is quasi-isometric to $G,$ $\adim G=\adim H.$ By Lemma~\ref{L:adimtree},$\adim H\leq 1.$ 
\end{proof}

The following theorem is a special case of \cite[Theorem 2.3]{DranishnikovSmith} in light of \cite[Theorem 2.1]{DranishnikovSmith}.

\begin{theorem}\label{T:adimsemi}
Let $G= H\ltimes K.$ If $G$ is finitely generable  and $K$ is locally finite, then $\adim G=\adim H$. 
\end{theorem}
\begin{proof}
As $G$ is finitely generable, so is $H$. Let $A$ be a finite generating set of $H$ and let $B$ be a finite subset of $K$ such that $X=A\cup B$ generates $G$. The inclusion $(H,\dist_A)\to(G,\dist_X)$ is an isometric embedding, and then $\adim (H,\dist_H)\leq \adim (G,\dist_X).$

We will denote $\dist_X$ and $\dist_A$ simply by $\dist.$ 

Let $h\in H$ and $k\in K.$ A geodesic expression for $hk$ has the form \begin{equation*}w_0,a_1,w_1,a_2,\dots, a_n,w_{n}\end{equation*} where $a_1,\dots,a_n$ are elements of $A^{\pm 1}$ and  $w_0,\dots,w_n$ are words in $B^{\pm 1}.$
As \begin{equation*}h=a_1\cdots a_n, \quad k=(a_1\cdots a_n)^{-1}w_0a_1w_1a_2\cdots a_nw_{n},\end{equation*}  hence  $\dist(1,hk)\geq \dist(1,h)$. Thus 
\begin{equation}\label{E:metric}
\text{for all } h\in H \text{ and all }k\in K, \dist(1,hk)\geq\dist(1,h).
\end{equation}

By \eqref{E:metric}, as $\dist$ is left $G$-invariant, we have that for $h,h'\in H$ and $k,k'\in K$
\begin{align*}\dist(hk,h'k') &= \dist(1,k^{-1}h^{-1}h'k')\\ &= \dist(1,(h^{-1}h')(k^{-1})^{h^{-1}h'}k')\\ &\stackrel{\eqref{E:metric}}{\geq} \dist(1,h^{-1}h')=\dist(h,h').
\end{align*}

Hence,
\begin{equation}\label{E:metric2}
\text{for all } h,h'\in H, \quad \dist(hK,h'K)= \dist(h,h').
\end{equation}

Fix $m>0$. We are going to show that $\adim (G,\dist) \leq \adim (H,\dist)$. By Observation~\ref{O:defadim}, there exist $r_H,z_H\colon H\to H$ such that  $h=z_H(h)r_H(h)$ for all $h\in H,$ $r_H(H)$ is finite and there is a $(1+\adim H)$-coloring of $z_H(H)$  such that for $z_H(h_1)\neq z_H(h_2)$ if $\col(z_H(h_1))=\col(z_H(h_2))$ then $\dist(h_1,h_2)>m.$ 

Let 
\begin{equation*}T=\{k\in K : \dist(h_1,h_2(k^{h_3}))\leq m  \text{ for some } h_1,h_2,h_3\in r_H(H)\}.\end{equation*} In particular $T$ is a finite subset of $K$, and hence $F=\gen{T}$ is a finite subgroup of $K$. Let $J$ be a transversal for the right $F$-multiplication action in $K.$ 

Let $i,j\in J.$ If $\dist(iFr_H(H),jFr_H(H))\leq m$ then there exist $f_i,f_j\in F,$ $h_i,h_j\in r_H(H)$ such that $\dist(if_ih_i,jf_jh_j)\leq  m.$ Hence \begin{equation*}\dist(h_i,f_i^{-1}i^{-1}jf_j h_j)=\dist(h_i,h_j(f_i^{-1}i^{-1}jf_j)^{h_j})\leq m.\end{equation*} By definition of $T$, $f_i^{-1}i^{-1}jf_j\in T\subset F$ and hence $i=j.$ This  proves
that
\begin{equation}\label{E:sep}
\text{for all } i,j\in J, \text{ if } i\neq j \text{ then } \dist(iFr_H(H),jFr_H(H))>m.
\end{equation}

For $g\in G$ there exist unique $h_g\in H,k_g\in K,j_g\in J$ and $f_g\in F$ such that $g=h_gk_g$ and $j_gf_g=\rconj{r_H(h_g)}{k_g},$ and hence 
\begin{align*}g=h_gk_g &= z_H(h_g)(r_H(h_g))k_g\\&= z_H(h_g)(r_H(h_g))k_g(r_H(h_g)^{-1})r_H(h_g)\\ &= z_H(h_g)j_gf_g(r_H(h_g)).\end{align*}

Thus the functions $z,r\colon G\to G,$ $z(g)=z_H(h_g)j_g$ and $r(g)=f_gr_H(h_g),$ are well defined and  $g=z(g)r(g)$ for all $g\in G$, $z(G)\subseteq Fr_H(H)$ is finite. Let $\col(z(g))=\col(z_H(h_g)).$ We have to show that for $g_1,g_2\in G$ if $z(g_1)\neq z(g_2)$ and $\col(z(g_1))=\col(z(g_2))$ then $\dist(g_1,g_2)>m.$

If $z_H(h_{g_1})\neq z_H(h_{g_2}),$ then $m<\dist(h_{g_1},h_{g_2})\leq \dist(g_1,g_2)$ by \eqref{E:metric2}.

If $z_H(h_{g_1})= z_H(h_{g_2}),$ then $j_{g_1}\neq j_{g_2}$ and 
\begin{align*}\dist(g_1,g_2)&\geq \dist(z_H(h_{g_1})j_{g_1}Fr_H(H),z_H(h_{g_2})j_{g_2}Fr_H(H))\\ &= \dist(j_{g_1}Fr_H(H),j_{g_2}Fr_H(H))\stackrel{\eqref{E:sep}}{\geq} m.\end{align*}
\end{proof}

\begin{definition}[Lamplighter group]\index{lamplighter group}
The {\it wreath product} $H\wr G$ of a group $H$ by a group $G$ is the semidirect product \begin{equation*}(\bigoplus_{g\in G} H)\rtimes G,\end{equation*} where $G$ acts  on $\oplus H$  by  multiplication on the coordinates.

Let $C_n$ denote the cyclic group of order $n$. The {\it lamplighter group } is the group $C_{2}\wr C_{\infty}.$ The standard presentation of the lamplighter group arises from the wreath product structure 
\begin{equation*} \langle a, t \mid a^2, [ t^m a t^{-m} , t^n a t^{-n} ], m, n \in \mathbb{Z} \rangle.\end{equation*}
\end{definition}

\begin{corollary}
The lamplighter group has asymptotic dimension 1.\qed
\end{corollary}

\begin{definition}
Let $\Gamma$ be a graph and $\gamma=(u_{[0\uparrow n]},e_{[1\uparrow n]})$ be a closed path in $\Gamma.$ 

A {\it diagram}  over $\Gamma$ for $\gamma$ is a pair $(D,\phi)$ where $D$ is a compact, planar, simply connected 2-complex  with a fixed embedding in $\R^2$ and $\phi$ is a graph map from the underlying graph of $D$ to $\Gamma$ such that the boundary of $D$ maps to $\gamma.$

A {\it triangular diagram} over $\Gamma$ for $\gamma$ is a diagram $(D,\phi)$ over $\Gamma$ for $\gamma$  where the boundary of each inner face of $D$ has at most 3 edges.

Let $\mathcal C=\prs{X}{R}$ be a presentation and $\Gamma$ be the Cayley graph of $\gp{X}{R}$ with respect to $X$. Let $w=x_{[0\uparrow n]}$ be a sequence in $X^{\pm 1}$ such that $\Pi x_{[0\uparrow n]}=1.$ Let $\mu$ be the closed path in $\Gamma$ given by $w.$

A {\it van Kampen diagram} over $\mathcal C$ for $w$ is a diagram $(D,\phi)$ over $\Gamma$ for $\mu$  such that for each inner face $d$ of $D,$ there exists a starting vertex and a direction for the closed path  bounding $d,$ say $(v_{[0\uparrow m]},f_{[1\uparrow m]})$, that satisfies that $(\phi(v_0)^{-1}\phi(v_1),\dots,\break \phi(v_{m-1})^{-1}\phi(v_m))$ is a freely reduced word of $R$.
\end{definition}
\begin{lemma}[van Kampen]
Let $\mathcal C=\prs{X}{R}$ be a presentation where each $r\in R$ is a freely reduced word. Let $w=x_{[0\uparrow n]}$ be a freely reduced word over $X$ with $\Pi x_{[0\uparrow n]}=1$ in $\gp{X}{R}.$
Then there exists a van Kampen diagram over $\mathcal C$ for $w.$ 
\qed
\end{lemma}
\begin{observation}\label{Obs:1seqtopathseq}
Let $v_{[0\uparrow n]}$ be an 1-sequence in $V\Gamma$. If $v_{[0\uparrow n]}$ is not a path sequence $\Gamma$, then there exists $i\in[0\uparrow n-1]$ such that $d(v_i, v_{i+1})=0.$ It follows that $v_i=v_{i+1}$ and, we can perform a $0$-reduction. Hence there exists a sequence of $0$-reductions of $v_{[0\uparrow n]}$ that produces a path sequence. 
\end{observation}
The main result of this section is that $G$ is virtually free if and only if $\adim(G)\leq 1$ and $G$ is finitely presentable. This result has been proved independently in \cite{FujiwaraWhyte}, \cite{Gentimis} and \cite{JanuszkiewiczSwiatkowski}. The proof given here is more elementary, and it relies on the results of the previous section.

\begin{theorem}\label{T:adim1}
Let $G$ be a finitely presentable group with $\adim G=1.$ Then there exists $M>0$ such that every $1$-sequence is $M$-triangulable.

\end{theorem}
\begin{proof}
Let $\prs{X}{R}$ be a finite presentation for $G$ and let $\Gamma$ be the Cayley graph of $G$ with respect to $X.$ 

If $R$ is empty, then $\Gamma$ is a tree and the result holds by Lemma~\ref{L:triangintree}. Hence we will assume that $R$ is non-empty and  every $r\in R$ is a freely reduced word.

Let $m=\max\{|r|:r\in R\}$ and let $z,r\colon G\to G$ such that $g=z(g)r(g)$ for all $g\in G$,  $r(G)$ is finite and there is a $2$-coloring for $z(G)$ such that for $g_1,g_2\in G$ if  $z(g_1)\neq z(g_2)$ and $\col(z(g_1))=\col(z(g_2))$ then $\dist(g_1,g_2)>m$.

Let $\Delta$ be the graph with vertex set $z(G)$ and an edge connecting two vertices $z(g),z(h)$ if and only if $z(g)\neq z(h)$ and $\dist(z(g)r(G),z(h)r(G))=1.$ It is straightforward to check that $z$ induces an (unoriented) graph quotient $\Gamma\to \Delta,$ denoted again by $z.$

Let $u_{[0\uparrow n]}$ be a $1$-sequence in $V\Gamma.$ If we show that the 1-sequence in $\Delta,$ $z(u_{[0\uparrow n]})$ is $1$-triangulable, then it follows that the sequence $u_{[0\uparrow n]}$ is $M$-triangulable with $M=(2\diam(r(G))+1).$ 

By Observation~\ref{Obs:1seqtopathseq}, we can assume that $u_{[0\uparrow n]}$ is a path sequence in $\Gamma.$ Let $\gamma=(u_{[0\uparrow n]},e_{[1\uparrow n]})$ be a reduced closed path in $\Gamma$. We may assume $n> 3$. By the van Kampen lemma, there exists a van Kampen diagram $(D_0,\phi_0)$ over $\prs{X}{R}$ for $\gamma.$ Hence $(D_0,z\circ \phi_0)$ is a diagram over $\Delta$ for $z(\gamma).$

We modify $D_0$ to obtain a triangular diagram over $\Delta$ for $z(\gamma).$

As each inner face $d$ of $D_0$ is labeled by some $r\in R,$ $\diam(\phi_0(\partial d))\leq m.$  We subdivide each inner face of $D_0$ into faces with at most 3 edges, following a triangulation of each $r\in R.$ Denote this new complex by $D$. 

For $x\in VD_0\cup ED_0$ let $\phi(x)=z(\phi_0(x)).$ For $e\in ED-ED_0,$ if $z(\iota e)=z(\tau e)$ let $\phi(e)=\phi(\iota e);$ if $\iota e\neq \tau e $ then $\dist(\iota e,\tau e)\leq m$ and hence $\dist(z^{-1}(\iota e),z^{-1}(\tau e))\leq 1,$ thus there exists a unique edge $f\in E\Delta$ with end-points $\{z(\iota e),z(\tau e)\}.$ Let $\phi(e)=f.$

Hence $(D, \phi)$ is a triangular diagram over $\Delta$ for $z(\gamma).$

\begin{figure}[ht]
\centerline{
\xymatrix{ D_0^{(0)}\cup D_0^{(1)}\ar[r]^{\phi_0}\ar[rd]^{i}&\Gamma\ar[r]^{z}&\Delta\\
& D^{(0)}\cup D^{(1)}\ar[ur]_{\phi}
}
}\caption{ The diagram $(D,z\circ\phi)$ over $\Delta$ for $\gamma$}\label{F:diag}
\end{figure}

We will show that there exists a triangular diagram $(D',\phi')$ over $\Delta$ for $z(\gamma)$ with no inner vertex. Note that if such diagram exists, then $z(\gamma)$ is $1$-triangulable in $\Delta$ and hence $\gamma$ is $M$-triangulable in $\Gamma.$

For $e\in ED,$ we call $e$ {\it strong} if $\col(\phi (\iota e))\neq\col(\phi(\tau e))$ otherwise we call it  {\it weak}. Equivalently, $e$ is strong if $\phi(e)\in E\Delta$ and weak if $\phi(e)\in V\Delta.$ 

Let $v\in VD,$ with $v$  not in the boundary of $D.$ 

\noindent\textbf{Construction 1:}
If there exists an edge $e\in ED$ such that $v\in \{\iota e,\tau e\}$ and  $e$ is a weak edge, then $\phi(e)$ is a vertex and $\phi$ is compatible with collapsing $e$ to a vertex in $D$. Collapsing the edge,  we obtain a new triangular diagram $(D',\phi')$ over $\Delta$ for $z(\gamma)$ with fewer inner vertices than $D.$

\noindent\textbf{Construction 2:}
Let $\Star{v}$ denote the set of $e\in ED$ such that $v\in \{\iota e,\tau e\}$. Suppose that all the elements of $\Star{v}$ are strong edges.

Since all the faces adjacent to $v$ have at most 3 edges, and if 2 edges of the face are strong, the other one should be weak, this implies that $z$ is constant on all the vertices of $\Star{v}$ different from $v$. Then we can modify $\phi$ to be constant in $\Star{v}$ as follows. Let $e\in \Star{v},$ $u\in \{\iota e,\tau e\}-\{v\}.$ Then we let $D'=D$ and $\phi'|_{D-\Star{v}\cup \{v\}}$ $=\phi$ and $\phi'(\Star{v}\cup \{v\})$ $=\phi(u).$  

Now applying the preceding construction, we reduce the number of inner vertices of $VD$.
\end{proof}

\begin{corollary}\label{Cor:adimvfg}
$G$ is virtually free if and only if $\adim G=1$ and $G$ finitely presentable.\qed
\end{corollary}

\begin{remark}
The Theorem~\ref{T:adim1} can be reformulated for graphs as follows: let $\Gamma$ be a graph with $\adim \Gamma \leq 1$ and $H_1(\Gamma)$ is generated by closed paths of length $<L$, then there exists $k>0$ such that every 1-sequence in $\Gamma$ is $k$-triangulable.
\end{remark}

\begin{corollary}\cite[Corollary 5]{BaumslagMiller}
Let $1\to K\to G\to H\to 1$ be a short exact sequence with $K$ locally finite, $G$ finitely generable and $H$ free. Then $G$ is finitely presentable if and only if $K$ is finite.

Moreover, if $K$ is finite then $G$ is virtually free.
\end{corollary}
\begin{proof}
As $H$ is free, there exists a section of $H\to G$ and  $G$ is a semidirect product.

By Theorem~\ref{T:adimsemi}, $\adim G=1.$

If $G$ is finitely presentable then, by Theorem~\ref{T:adim1}, $G$ is virtually free. Hence, there exists a short exact sequence $1\to N\to G\stackrel{f}{\to} \Phi\to 1$ with $N$ free and $\Phi$ finite. As $K\cap N=\{1\},$ the restriction of $f$ to $K$ has trivial kernel, and hence $K$ is finite.

If $K$ is finite, then as $G$ is a semidirect product of finitely presentable groups, $G$ is finitely presentable. Now, by Theorem~\ref{T:adim1}, $G$ is virtually free.
\end{proof}

\noindent{\textbf{\Large{Acknowledgments}}}

The author is grateful to Warren Dicks, Pere Menal, Sarah Rees and Nickolas Wright for several helpful conversations.
\medskip

\footnotesize

The research of the  author was
jointly funded by the MEC (Spain) and the EFRD (EU)
through Projects MTM2006-13544 and MTM2008-01550.

\bibliographystyle{amsplain}

\bibliographystyle{amsplain}

\textsc{School of Mathematics,
University of Southampton,
Southampton
SO17 1BJ
}

\emph{E-mail address}{:\;\;}\texttt{Y.Antolin-Pichel@soton.ac.uk}

\end{document}